\documentclass{amsart}
\usepackage{amsmath, amsthm, amssymb, amscd, epsfig}
\usepackage{pinlabel}

\begin{document}

\newtheorem{theorem}{Theorem}[section]
\newtheorem{lemma}[theorem]{Lemma}
\newtheorem{proposition}[theorem]{Proposition}
\newtheorem{corollary}[theorem]{Corollary}
\newtheorem{conjecture}[theorem]{Conjecture}
\newtheorem{question}[theorem]{Question}
\newtheorem{problem}[theorem]{Problem}
\newtheorem*{claim}{Claim}
\newtheorem*{criterion}{Criterion}
\newtheorem*{stability_thm}{Stability Theorem~\ref{stability_theorem}}

\theoremstyle{definition}
\newtheorem{definition}[theorem]{Definition}
\newtheorem{construction}[theorem]{Construction}
\newtheorem{notation}[theorem]{Notation}

\theoremstyle{remark}
\newtheorem{remark}[theorem]{Remark}
\newtheorem{example}[theorem]{Example}

\numberwithin{equation}{subsection}

\def\Z{\mathbb Z}
\def\R{\mathbb R}
\def\Q{\mathbb Q}
\def\H{\mathbb H}
\def\E{\mathcal E}
\def\M{\mathcal M}

\def\freq{\textnormal{freq}}
\def\cl{\textnormal{cl}}
\def\scl{\textnormal{scl}}
\def\homeo{\textnormal{Homeo}}
\def\rot{\textnormal{rot}}
\def\area{\textnormal{area}}

\def\Id{\textnormal{Id}}
\def\PSL{\textnormal{PSL}}
\def\til{\widetilde}

\title{Immersed surfaces in the modular orbifold}
\author{Danny Calegari}
\address{Department of Mathematics \\ Caltech \\
Pasadena CA, 91125}
\email{dannyc@its.caltech.edu}
\author{Joel Louwsma}
\address{Department of Mathematics \\ Caltech \\
Pasadena CA, 91125}
\email{louwsma@caltech.edu}
\date{\today}

\begin{abstract}
A hyperbolic conjugacy class in the modular group $\PSL(2,\Z)$ corresponds 
to a closed geodesic in the modular orbifold. Some of these geodesics 
virtually bound immersed surfaces, and some do not; the distinction is 
related to the polyhedral structure in the unit ball of the stable 
commutator length norm. We prove the following stability theorem: for every 
hyperbolic element of the modular group, the product of this element with 
a sufficiently large power of a parabolic element is represented by a geodesic 
that virtually bounds an immersed surface. 
\end{abstract}

\maketitle

\section{Introduction}

In many areas of geometry, it is important to understand which
immersed curves on a surface bound immersed subsurfaces. Such questions arise (for example) in topology,
complex analysis, contact geometry and string theory. In \cite{Calegari_faces} it was shown
that studying isometric immersions between hyperbolic surfaces with geodesic boundary gives
insight into the polyhedral structure of the second bounded cohomology of a free group
(through its connection to the {\em stable commutator length} norm, defined in \S~\ref{scl_subsection}).

If $\Sigma$ is a (complete) noncompact oriented hyperbolic orbifold, a (hyperbolic)
conjugacy class $g$ in $\pi_1(\Sigma)$ is represented by
a unique geodesic $\gamma$ on $\Sigma$. The {\em immersion problem} asks when there is an oriented surface $S$ and an
orientation-preserving immersion $i:S \to \Sigma$ taking $\partial S$ to $\gamma$ in an orientation-preserving
way. If the problem has a positive solution one
says that $\gamma$ {\em bounds an immersed surface}. The immersion problem is complicated 
by the fact that there are examples of curves $\gamma$ that do not bound immersed surfaces, but have
finite (possibly disconnected) covers that do bound; i.e. there is an immersion $i:S \to \Sigma$ as above
for which $\partial S \to \gamma$ factors through a {\em covering map} of some positive degree. In
this case we say $\gamma$ {\em virtually bounds} an immersed surface.

One can extend the virtual immersion problem in a natural way to finite rational formal sums of geodesics representing
zero in (rational) homology. Formally, one defines the real vector space $B_1^H(G)$ of {\em homogenized $1$-boundaries},
(see \S~\ref{scl_subsection}), where $G=\pi_1(\Sigma)$. It is shown in \cite{Calegari_faces} that the set 
of rational chains in $B_1^H$ for which the virtual immersion problem has a positive solution are
precisely the set of rational points in a closed convex rational polyhedral cone with nonempty interior. 
This fact is connected in a deep way with Thurston's characterization (see \cite{Thurston}) of the
set of classes in $H_2(M;\Z)$ represented by fibers of fibrations of a given $3$-manifold $M$ over $S^1$. 
It can be used to give a 
new proof of symplectic rigidity theorems of Burger-Iozzi-Wienhard \cite{Burger_Iozzi_Wienhard}, 
and has been used by Wilton \cite{Wilton} in his work on Casson's conjecture.
Evidently understanding the structure of the set of solutions of the virtual immersion
problem is important, with potential applications to many areas of mathematics. Unfortunately,
this set is apparently very complicated, even for very simple surfaces $\Sigma$.

\medskip

The purpose of this paper is to prove the following Stability Theorem:

\begin{stability_thm}
Let $v$ be any hyperbolic conjugacy class in $\PSL(2,\Z)$, represented by a string $X$ of
positive $R$'s and $L$'s. Then for all sufficiently large $n$,
the geodesic in the modular orbifold $\M$ corresponding to the stabilization $R^n X$ 
virtually bounds an immersed surface in $\M$.
\end{stability_thm}

This theorem proves the natural analogue of Conjecture~3.16 
from \cite{Calegari_faces}, with $\PSL(2,\Z)$ in place of the free group $F_2$.

It follows from the main theorems of \cite{Calegari_faces}
that the elements $v_n \in \PSL(2,\Z)$ corresponding to the stabilizations of $v$
as above satisfy $\scl(v_n) = \rot(v_n)/2$,
where $\rot$ is the {\em rotation quasimorphism} on $\PSL(2,\Z)$, and $\scl$ denotes the
{\em stable commutator length} (see \S~\ref{background_section} for details). Under the natural central
extension $\phi: B_3 \to \PSL(2,\Z)$ where $B_3$ denotes the $3$ strand braid group,
there is an equality $\scl(b) = \scl(\phi(b))$ for all $b \in [B_3,B_3]$;
consequently we derive an analogous stability theorem for stable commutator length in $B_3$.

\medskip

We give the necessary background and motivation in \S~\ref{background_section}.
Theorem~\ref{stability_theorem} is proved in \S~\ref{proof_section}. In \S~\ref{generalize_section}
we generalize our main theorem to $(2,p,\infty)$-orbifolds for any $p\ge 3$, and discuss
some related combinatorial problems and a connection to a problem in theoretical computer science. 
Finally, in \S~\ref{experimental_section} we describe the results of some computer experiments.

\section{Background}\label{background_section}

We recall here some standard definitions and facts for the convenience of the reader.

\subsection{The modular group}
The modular group $\PSL(2,\Z)$ acts discretely and with finite covolume on $\H^2$, and the quotient
is the {\em modular orbifold} $\M$, which can be thought of as a triangle orbifold of type 
$(2,3,\infty)$.

Every element of $\PSL(2,\Z)$ either has finite order, or is parabolic, or is conjugate to a product of the
form $R^{a_1}L^{b_1}R^{a_2} \cdots L^{a_m}$ where the $a_i$ and $b_i$ are positive integers, and
where $L$ and $R$ are represented by the matrices
$$L = \begin{pmatrix} 1 & 0 \\ 1 & 1 \\ \end{pmatrix} \quad R = \begin{pmatrix} 1 & 1 \\ 0 & 1 \\ \end{pmatrix}$$
(the parabolic elements are conjugate into the form $R^a$ or $L^b$).

The group $\PSL(2,\Z)$ is abstractly isomorphic to the free product $\Z/2\Z * \Z/3\Z$. 

\subsection{Stable commutator length}\label{scl_subsection}

For a basic introduction to stable commutator length, see \cite{Calegari_scl}, especially Chapters 2 and 4.

If $G$ is a group, and $g \in [G,G]$, the {\em commutator length} of $g$ (denoted $\cl(g)$) is the
smallest number of commutators in $G$ whose product is $g$, and the {\em stable commutator length}
$\scl(g)$ is the limit $\scl(g):=\lim_{n \to \infty} \cl(g^n)/n$.

Stable commutator length extends in a natural way to a function on $B_1(G)$, the space of
real group $1$-boundaries (i.e. real group $1$-chains representing $0$ in homology with real
coefficients) and descends to a pseudo-norm (which can be thought of as a kind of
relative Gromov-Thurston norm) on the quotient $B_1^H(G):=B_1(G)/\langle g-hgh^{-1}, g^n-ng\rangle$.

The dual of this space (up to scaling by a factor of $2$) is $Q(G)/H^1(G)$ --- 
i.e. homogeneous quasimorphisms on $G$ modulo homomorphisms, with the defect norm. This duality
theorem --- known as {\em Generalized Bavard Duality} --- is proved in \cite{Calegari_scl}, \S~2.6.
A special case of this theorem, proved by Bavard in \cite{Bavard}, says that for any $g \in [G,G]$
there is an equality $\scl(g) =\sup_\phi \phi(g)/2$
where the supremum is taken over all homogeneous quasimorphisms $\phi \in Q(G)$ normalized to have
defect $D(\phi)$ equal to $1$. A quasimorphism $\phi$ (with $D(\phi)=1$)
is {\em extremal} for $g$ if $\scl(g)=\phi(g)/2$. 

The theory of stable commutator length has deep connections to dynamics, group theory, geometry,
and topology; however, although in principle this function contains a great deal of information,
it is notoriously difficult to extract this information, and to interpret it geometrically.
It is a fundamental question in any group $G$
to calculate $\scl$ on chains in $B_1^H(G)$, and to determine extremal quasimorphisms for such
elements. Conversely, given a homogeneous quasimorphism $\phi \in Q(G)$ (especially one with
some geometric ``meaning''), it is a fundamental
question to determine the (possibly empty) cone in $B_1^H(G)$ on which $\phi$ is extremal.

\subsection{Rotation quasimorphism}

If $G$ is a word-hyperbolic group, $\scl$ is a genuine norm on $B_1^H(G)$. Moreover, it is shown
in \cite{Calegari_rational} and \cite{Calegari_faces} that if $G$ is
virtually free, the unit ball in this norm is a rational polyhedron, and there are codimension one
faces associated to realizations of $G$ as the fundamental group of a complete oriented hyperbolic
orbifold.

Dual to each such codimension one face is a unique extremal vertex of the unit ball in $Q/H^1$. In
our case, $\PSL(2,\Z)$ may be naturally identified with the fundamental group of the modular orbifold.
The unique homogeneous quasimorphism dual to this realization, scaled to have defect $1$, is the
{\em rotation quasimorphism}, denoted $\rot$. This rotation function is very closely related to
the {\em Rademacher $\varphi$ function}, which arises in connection with Dedekind's $\eta$ function, and
is studied by many authors, e.g. \cite{Atiyah, Ghys, Kirby_Melvin, Rademacher, Hirzebruch_Zagier}, and so on
(see especially \cite{Ghys}, \S~3.2 for a discussion most closely connected to the point of
view of this paper). In fact, up to a constant, the rotation quasimorphism is the {\em homogenization} of the
Rademacher function; i.e. $\rot(g) = \lim_{n \to \infty} \varphi(g^n)/6n$.

\begin{figure}[htpb]
\labellist
\small\hair 2pt
%\pinlabel $0$ at 310 0
\endlabellist
\centering
\includegraphics[scale=1]{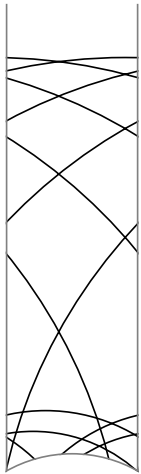}
\caption{$\gamma$ cuts a fundamental domain into regions $R_i$}\label{R7L2RL_jigsaw}
\end{figure}

The simplest way to define this function (at least on hyperbolic elements of $\PSL(2,\Z)$)
is as follows. Associated to a hyperbolic
conjugacy class $g \in \PSL(2,\Z)$ is a geodesic $\gamma$ on $\M$. The geodesic $\gamma$ cuts $\M$
up into complementary regions $R_i$ (see Figure~\ref{R7L2RL_jigsaw} for an example). 
Join each region $R_i$ to the cusp by a proper ray $\alpha_i$,
and define $n_i$ to be the {\em signed} intersection number $n_i=\alpha_i \cap \gamma$. Then
$\rot(g)=\frac 1 {2\pi} \sum_i n_i\, \area(R_i)$. In other words, up to a factor of $2\pi$, the
rotation number is the {\em algebraic area} enclosed by $\gamma$. Algebraically, if the conjugacy
class of $g$ has a factorization of the form $R^{a_1}L^{b_1}\cdots L^{b_n}$ then 
$\rot(g) = \frac 1 6 (\sum a_i - \sum b_i)$; see \cite{Kirby_Melvin} 1.5--6, or \cite{Rademacher} equation 70
(Rademacher denotes $\varphi$ and $6\,\rot$ by $\Phi$ and $\Psi$ respectively).

By Bavard duality, one has $\scl(g) \ge \rot(g)/2$ for every $g$. Moreover, it is shown in 
\cite{Calegari_faces} (for arbitrary free groups, though the proof easily generalizes to virtually
free groups) that equality is achieved if and only if the 
geodesic representative $\gamma$ of $g$ virtually bounds an immersed surface. That is, if and only 
if there is a hyperbolic surface $S$ and an isometric immersion $S \to \M$ wrapping 
$\partial S$ some (positive) number $n$ times around $\gamma$. Topologically, one can think of
the problem of constructing such an immersed surface as a kind of {\em jigsaw puzzle}: 
one takes $n\cdot n_i$ copies of each region $R_i$, where $n$, $n_i$ and $R_i$ are as above, 
and tries to glue them up compatibly with their tautological embeddings in $\M$ in such a way 
as to produce a smooth orbifold with geodesic boundary. Evidently, a necessary condition is that
the $n_i$ are all non-negative. However, this necessary condition is not sufficient.

\medskip

In \cite{Calegari_faces} it was observed experimentally that for many words $w \in [F_2,F_2]$,
geodesics on a hyperbolic once-punctured torus corresponding to conjugacy classes of the form
$[a,b]^nw$ all virtually bound immersed surfaces for sufficiently large $n$, and it was
conjectured (Conjecture~3.16) that this holds in general. Our main theorem (Theorem~\ref{stability_theorem}
below) proves the natural analogue of this conjecture with the free group $F_2$ replaced by the
virtually free group $\PSL(2,\Z)$.

\subsection{Braid group}

The braid group $B_3$ is a central extension of $\PSL(2,\Z)$. Under this projection, the
standard braid generators $\sigma_1$ and $\sigma_2$ are taken to $R^{-1}$ and $L$ respectively.
It is straightforward to show that for any $b \in [B_3,B_3]$ the (stable) commutator length of $b$
is equal to the (stable) commutator length of its image in $\PSL(2,\Z)$. Consequently, our
main theorem shows that the rotation quasimorphism is extremal for a sufficiently large 
stabilization of any element of $[B_3,B_3]$.

\section{Proof of theorem}\label{proof_section}

The purpose of this section is to prove the following theorem:

\begin{theorem}[Stability theorem]\label{stability_theorem}
Let $v$ be any hyperbolic conjugacy class in $\PSL(2,\Z)$, represented by a string $X$ of
positive $R$'s and $L$'s. Then for all sufficiently large $n$,
the geodesic in the modular orbifold $\M$ corresponding to the stabilization $R^n X$ 
virtually bounds an immersed surface in $\M$.
\end{theorem}

The proof will occupy the remainder of the section.

\medskip

The conjugacy class of $v$ has a representative of the form $R^{a_1}L^{b_1}R^{a_2}\cdots L^{b_n}$
where the $a_i,b_i$ are all positive integers.
We will show that the geodesic $\gamma$ corresponding to $v$ virtually bounds an immersed surface
providing $a_1$ is sufficiently big compared to $\sum_{i\ne 1} a_i + \sum_i b_i$. Evidently the theorem follows from this.
We fix the notation $N=a_1$ and $N' = \sum_{i\ne 1} a_i + \sum_i b_i$ in the sequel, and we prove
the theorem under the hypothesis $N \ge 3N' + 11n + 3$ (note that there is no suggestion that this
inequality is sharp).

\begin{figure}[htpb]
\labellist
\small\hair 2pt
%\pinlabel $\nu$ at -20 75
\endlabellist
\centering
\includegraphics[scale=0.5]{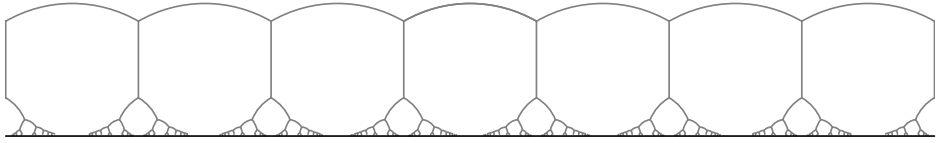}
\caption{The tree $\til{\sigma}$ in the upper half-space model} \label{tildesigma}
\end{figure}

In the modular orbifold $\M$, let $\sigma$ denote the embedded geodesic segment running between the
orbifold points of orders $2$ and $3$. The preimage of $\sigma$ in the universal cover $\H^2$ is a
regular $3$-valent tree, which we denote $\til{\sigma}$; see Figure~\ref{tildesigma}.

In the upper half-space model of $\H^2$, let $W$ be the closure of the
complementary component of $\H^2 - \til{\sigma}$ stabilized by the translation $z \to z+1$. 
Then $\partial W$ is a collection of circular arcs, whose vertices are the complex numbers
$e^{2\pi i/3}+n$ for $n \in \Z$. We call these arcs the {\em segments} of $\partial W$. Every
segment of $\partial W$ (orbifold) double covers the interval $\sigma$ in $\M$.
In the sequel we use the abbreviation $\omega = e^{2\pi i/3}$.

\subsection{Arcs and subwords}

The arc $\sigma$ cuts $\gamma$ into a collection of geodesic segments which 
correspond approximately to the $R^{a_i}$ and $L^{b_i}$ terms in the expression of $v$, in the
following way. 

After choosing a base point and an orientation, a word $v$ in the $L$'s and $R$'s 
determines a simplicial path in $\til{\sigma}$, where $L$ indicates
a ``left turn'', and $R$ indicates a ``right turn'', reading the word from right to left. 
The letters $R$ or $L$ correspond to the
vertices of this path, and a string of the form $RL^bR$ (resp. $LR^aL$) corresponds to a segment of length
$b+1$ (resp. $a+1$) which, after translation by some element of $\PSL(2,\Z)$, we can arrange to be contained
in $\partial W$ as a string of $b+1$ consecutive arcs moving to the right (resp. $a+1$ consecutive
arcs moving to the left).

The bi-infinite power $\dot{v}:=\cdots vvv\cdots$
determines a bi-infinite path $P(v)$ in $\til{\sigma}$, which is a quasigeodesic in $\H^2$ a bounded distance
from the geodesic representative of an axis of (some conjugate of) $v$. We may thus crudely associate lifts of
segments $\gamma_j$ to $W$ with such subwords.

If we translate $P(v)$ so that the segment corresponding to $L^b$ starts at the
vertex $\omega$ of $\til{\sigma}$, then this segment ends at $\omega + b+1$. Moreover,
the endpoints of $P(v)$ on the real axis are contained in the intervals $(-1,0)$ and $(b,b+1)$.
Let $\til{\gamma}$ be the infinite geodesic with the same endpoints as $P(v)$. Then the
intersection of $\til{\gamma}$ with $W$ is either empty (for which $b=1$ is necessary but
not sufficient) or consists of two points, one in the segment of
$\til{\sigma}$ from $\omega$ to $\omega\pm 1$, the other in the segment of $\til{\sigma}$ from
$\omega+b$ to $\omega+b\pm 1$ where the $\pm 1$ depends in each case on the rest of the word $v$
(the degenerate case that $\til{\gamma}$ passes through one or two vertices of $\til{\sigma}$ is
allowed).

In our case of interest, this intersection $\til{\gamma} \cap W$ projects to 
the segment $\beta_i$ of $\gamma$ corresponding to an $L^{b_i}$ subword. Similarly, $\alpha_i$
segments of $\gamma$ correspond to $R^{a_i}$ subwords, with the caveat that
some $L^{b_i}$ or $R^{a_i}$ subwords with $a_i$ or $b_i=1$ may not correspond to a segment 
of $\gamma$ at all.

\begin{example}[$R^7L^2RL$ part 1]
We will illustrate the main points of our construction in a particular case. 

\begin{figure}[htpb]
\labellist
\small\hair 2pt
\pinlabel $P(v)$ at 500 130
\pinlabel $\longleftarrow$ at 500 75
\endlabellist
\centering
\includegraphics[scale=0.25]{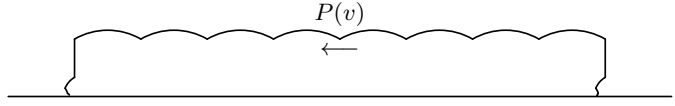}
\caption{A bi-infinite path $P(v)$ associated to $v=R^7L^2RL$. At each vertex, $P$ turns
left or right according to the letters of $\dot{v}$ (read from right to left). The path
$P(v)$ is chosen so that a segment corresponding to $R^7$ is contained in $\partial W$.}\label{R7L2RL_P}
\end{figure}

Let $v$ correspond to the conjugacy class $R^7L^2RL$ which
satisfies $\scl(v)=5/12$ and $\rot(v)=5/6$. A matrix representative in $\PSL(2,\Z)$ for $v$ is
$\left(\begin{smallmatrix} 37 & 22 \\ 5 & 3 \\ \end{smallmatrix}\right)$.
A bi-infinite path $P(v)$ is illustrated in Figure~\ref{R7L2RL_P} and the corresponding axis
$\til{\gamma}$ in Figure~\ref{R7L2RL_gamma}.

\begin{figure}[htpb]
\labellist
\small\hair 2pt
\pinlabel $R^7$ at 350 220
\pinlabel $\longleftarrow$ at 348 195
\pinlabel $\text{prefix:}$ at 80 -14
\pinlabel $\dot{L}$ at 175 -10
\pinlabel $\dot{R}L$ at 125 -10
\pinlabel $\text{suffix:}$ at 490 -12
\pinlabel $\dot{L}$ at 525 -10
\pinlabel $L\dot{R}$ at 575 -10
\endlabellist
\centering
\includegraphics[scale=0.5]{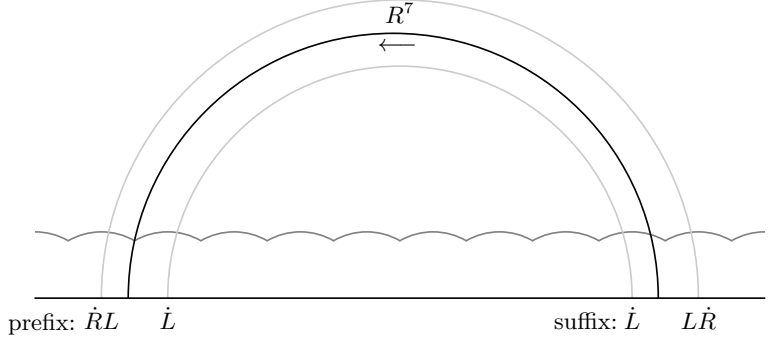}
\caption{An axis $\til{\gamma}$ obtained by straightening $P(v)$.
The ``barrier'' circles are associated to the bi-infinite words $\dot{R}LR^7L\dot{R}$
and $\dot{L}R^7\dot{L}$. The part of $\til{\gamma}$ contained in $W$ is a 
lift of $\alpha_1$.}\label{R7L2RL_gamma}
\end{figure}

The geodesic $\gamma$ is cut by $\sigma$ into three segments $\alpha_1$, $\beta_1$, $\alpha_2$, 
corresponding to the subwords $R^7$, $L^2$ and $R$.
\end{example}

\subsection{Lifts and surfaces}\label{lifts_surfaces_subsection}
For each $\alpha_i$ or $\beta_i$ we choose lifts $\til{\alpha}_i$ and $\til{\beta}_i$ properly
contained in $W$ subject to the following {\em lifting conditions}:
\begin{enumerate}
\item{the lifts are disjoint, and no two lifts intersect the same segment of $\partial W$}
\item{there are exactly five consecutive segments of $\partial W$ between
consecutive $\til{\beta}_i$}
\item{the $\til{\alpha}_i$ and $\til{\beta}_i$ are not nested (i.e. they cobound disjoint disks with $\partial W$)
except that the $\til{\beta}_i$ are all contained ``under'' the lift $\til{\alpha}_1$, so that the
leftmost vertex of $\til{\alpha}_1$ and the leftmost vertex of $\til{\beta}_1$ intersect 
segments of $\partial W$ separated by exactly five other segments of $\partial W$.}
\end{enumerate}

For each $i \ne 1$ let $S_i$ be the subsurface of $W$ bounded by $\til{\alpha}_i$ and $\partial W$. Let
$T$ be the subsurface of $W$ ``above'' all the $\til{\beta}_i$ and ``below'' $\til{\alpha}_1$.
We will build our immersed surface from the $S_i$ and $T$, glued up suitably along their
intersection with $\partial W$. Let $S$ denote the (disjoint) union $S = \cup_i S_i \cup T$.
The boundary of $S$ comes in two parts: the part of
the boundary along the $\til{\alpha}_i$ and $\til{\beta}_i$ (we denote this part of the boundary
$\partial_\gamma S$) and the part along $\partial W$ (we denote this part by $\partial_W S$).
Furthermore, $\partial_W S$ decomposes naturally into segments which are the intersection with
the segments of $\partial W$. There are two kinds of such segments: {\em entire} segments (those
corresponding to an entire segment of $\partial W$ contained in $\partial S$), and
{\em partial} segments (those corresponding to a segment of $\partial W$ containing an endpoint of
some $\til{\alpha}_i$ or $\til{\beta}_i$ in its interior). We also allow the case of a degenerate
partial segment, consisting of a single vertex of $\partial W \cap \til{\alpha}_i$ or
$\partial W \cap \til{\beta}_i$.

By construction, the partial segments of $\partial_W S$ come in oppositely oriented pairs, ending
on pairs of points of $\partial_\gamma S \cap \partial_W S$ projecting to the same point in $\gamma$.
We glue up such pairs of partial segments, producing a surface $S'$. Under the covering projection $\H^2 \to \M$,
the surface $S$ immerses in $\M$ in such a way that the immersion extends to an immersion of $S'$.
The $\partial_\gamma$ components glue up to produce a smooth boundary component $\partial_\gamma S'$ 
which wraps once around $\gamma$ in $\M$. The other part of $\partial S'$, which by abuse of
notation we denote $\partial_W S'$, is a union of connected components, each of which is tiled by
entire segments of $\partial W$. At each vertex corresponding to an end vertex of a partial segment,
the segments of $\partial_W S'$ meet at an angle of $4\pi/3$. At every other vertex the segments
meet at an angle of $2\pi/3$. We say such vertices are of {\em type $2$} and {\em type $1$} respectively.

To complete the construction of an immersed surface virtually bounding $\gamma$ (and thereby
completing the proof of Theorem~\ref{stability_theorem}) we must show how to glue up $S'$ 
by identifying segments of $\partial_W$ to produce a smooth surface. Such a
surface immerses in $\M$, and is extremal for $\gamma$. In fact, technically it is easier to glue up
$S'$ to produce a smooth {\em orbifold}, containing orbifold points of order $2$ and $3$ that map
to the corresponding orbifold points in $\M$. Such an orbifold is finitely covered by a smooth surface
virtually bounding $\gamma$.

\begin{example}[$R^7L^2RL$ part 2]
With notation as in part 1, we choose lifts 
$\til{\alpha}_1$, $\til{\beta}_1$ and $\til{\alpha}_2$ as indicated in Figure~\ref{R7L2RL_lifts}.

\begin{figure}[htpb]
\labellist
\small\hair 2pt
\pinlabel $\til{\alpha}_1$ at 320 220
\pinlabel $\til{\alpha}_2$ at 640 80
\pinlabel $\til{\beta}_1$ at 390 87
\pinlabel $a_1$ at 125 40
\pinlabel $c_1$ at 437 40
\pinlabel $b_2$ at 350 43
\pinlabel $b_1$ at 604 43
\pinlabel $a_2$ at 678 40
\pinlabel $c_2$ at 513 40
\endlabellist
\centering
\includegraphics[scale=0.5]{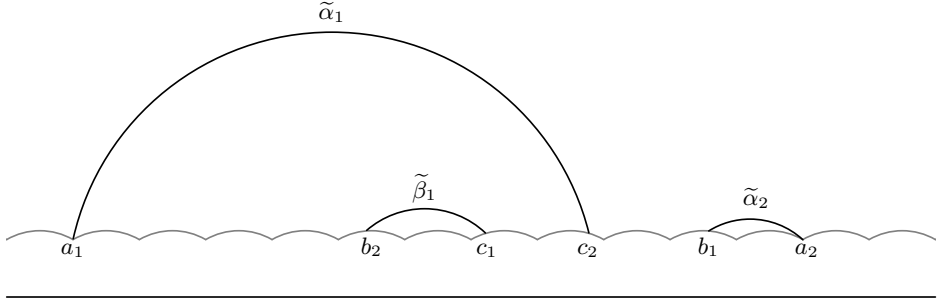}
\caption{Lifts $\til{\alpha}_1$, $\til{\beta}_1$ and $\til{\alpha}_2$. The region $T$ is
between $\til{\alpha}_1$ and $\til{\beta}_1$, and the region $S_2$ is under $\til{\alpha}_2$.
The endpoints of the $\til{\alpha}_i$ and $\til{\beta}_1$ are glued up by gluing
$a_1$ to $a_2$, $b_1$ to $b_2$, and $c_1$ to $c_2$.}\label{R7L2RL_lifts}
\end{figure}

Notice that the exponent $7$ is too small: there is not enough room for $\til{\beta}_1$
under $\til{\alpha}_1$ without putting two endpoints ($c_1$ and $c_2$) on adjacent partial
segments. The result of gluing up these adjacent segments produces an orbifold point of
order $3$ in the interior of $S'$.

Furthermore, there are a pair of ``degenerate'' partial segments, consisting
of the points $a_1$ and $a_2$. Identifying these points produces a non-manifold point in
$S'$, where the two components $\partial_\gamma S'$ and $\partial_W S'$ meet, at a smooth
point on $\partial_\gamma S'$, and at a point with angle $4\pi/3$ (i.e. a point of type $2$)
on $\partial_W S'$. This non-manifold point will become an ordinary manifold point after 
we glue up $\partial_W S'$.

The $a_i$ and the endpoints of the partial segments containing the $b_i$ 
glue up into two adjacent type $2$ points on $\partial_W S'$, and the remaining
three vertices of $\partial_W S$ give rise to three adjacent type $1$ points on $\partial_W S'$.
Thus the vertices on $\partial_W S'$ are of type $11221$ in cyclic order.

The interior of the $12$ segment can be folded up, creating an orbifold point of order $2$,
and the other four segments identified in pairs (pairing $11$ with $22$), 
creating another orbifold point of order $3$ corresponding to the ``unpaired'' $1$ vertex. 

The result is a smooth orbifold $S''$ with three orbifold points of orders $2,3,3$, which immerses
in $\M$ with boundary wrapping once around $\gamma$. A finite (orbifold) cover of $S''$ is a
genuine surface, which $\gamma$ virtually bounds.
\end{example}

\subsection{Combinatorics}

We have seen in general that $\partial_W S'$ is determined by combinatorial 
data consisting of a finite collection of circularly ordered sequences of 
$1$'s and $2$'s, which we call {\em circles}. We write such a circle
as an ordered list of $1$'s and $2$'s, where two such ordered lists define the same
circle (denoted $\sim$) if they differ by a cyclic permutation. Hence $211\sim 121 \sim 112$
and so on. A consecutive string of $1$'s and $2$'s contained in a circle is bracketed by a
dot on either side; hence $\cdot 12\cdot$ is a sequence in the circle $122$ and so on.

The $2$'s correspond to vertices of $\partial_W S$
which are also vertices of $\partial W$, and are contained in partial edges; whereas the $1$'s 
correspond to the other vertices of $\partial_W S$ which are also vertices of $\partial W$. Hence
the total number of $2$'s is equal to the number of segments of $\gamma - \sigma$, which is at most $n$.
The total number of $1$'s is at most $N+N'+n$. 

The only combinatorial properties of the circles we need are the following:
\begin{enumerate}
\item{each circle contains at least one $2$, and consequently there are at most $n$ circles
(this is immediate from the construction);}
\item{some circle contains a sequence of at least $N-N'-7n$ consecutive $1$'s, where $N$ is
large compared to $N'$ and $n$
(this follows from lifting conditions (2) and (3)); and}
\item{each circle contains a string of at least two consecutive $1$'s (this follows from lifting
condition (2)).}
\end{enumerate}

We refer to the sequence of at least $N-N'-7n$ consecutive $1$'s informally as 
{\em the big $1$ sequence}. Providing $N$ is sufficiently big compared to $N'$ and $n$ ---
equivalently, providing the big $1$ sequence is sufficiently long --- we can completely glue up 
$\partial_W S'$ as an orbifold. We now explain how to do this.

\medskip

The argument consists of a sequence of reductions to simpler and simpler combinatorial configurations.
These reductions are described using a pictorial calculus whose meaning should be self-evident.

\begin{figure}[htpb]
\labellist
\small\hair 2pt
\pinlabel $\xrightarrow{1\text{-handle}}$ at 75 33
\pinlabel $1$ at 13 21
\pinlabel $1$ at 13 41
\pinlabel $1$ at 47 21
\pinlabel $1$ at 47 41
\pinlabel $2$ at 115 15
\pinlabel $2$ at 115 45
\endlabellist
\centering
\includegraphics[scale=1]{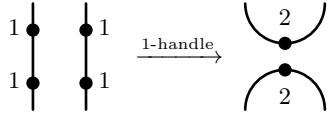}
\caption{The $1$-handle move}\label{1handle}
\end{figure}

The first reduction consists of taking a pair of $11$'s and identifying the segments
they bound. If the $11$'s are on different circles, these two circles become amalgamated
into a single circle. We call this {\em the $1$-handle move}; see Figure~\ref{1handle}.

Hence, by bullets (2) and (3), by
applying the $1$-handle move at most $n$ times, using up at most $2n$ of the $1$'s in the big $1$ 
sequence in the process, we can reduce to the case of a single circle. This circle has the
form $1^m\nu$ where $1^m$ is the big $1$ sequence, and $\nu$ stands for some sequence of $1$'s and
$2$'s. Since $N\ge 3N' + 11n + 3$, the big $1$ sequence in $\partial S'$ has length at least $2N' + 6n + 3$.
Using up at most $2n$ of these $1$'s gives $m \ge 2N' + 4n + 3$. On the other hand, the length
of $\nu$ is at most $N'+2n$ by construction.

We introduce three other simple moves:
\begin{enumerate}
\item{a single $\cdot 11\cdot$ is folded up, producing an interior orbifold point of order $2$ and a single $2$ vertex;
i.e. $\cdot 11\cdot \to \cdot 2\cdot$; or}
\item{the two adjacent $\cdot 11\cdot$ segments in a $\cdot 111\cdot$ are identified, producing an interior orbifold point of order
$3$ and a single $2$ vertex; i.e. $\cdot 111\cdot \to \cdot 2\cdot$; or}
\item{the middle $\cdot 12\cdot$ segment of a $\cdot 1121\cdot$ is folded up, producing an interior orbifold point of order $2$ and
a single $2$ vertex; i.e. $\cdot 1121\cdot \to \cdot 2\cdot$ (also $\cdot 1211\cdot \to \cdot 2\cdot$).}
\end{enumerate}
The first two moves are special cases of the $1$-handle move, where the two $\cdot 11\cdot$ segments being glued
are not disjoint. By means of repeated applications of the $\cdot 11\cdot\to \cdot 2\cdot$ move, a string of 
at most $2k$ consecutive $1$'s can be reduced to any string of length $k$.

Associated to any sequence of $1$'s and $2$'s is its {\em complement}, obtained by reversing the order
of the sequence and replacing each $1$ by a $2$ and conversely. We denote the complement of a
sequence $\nu$ by $\nu^c$. We transform a subset of the
big $1$ sequence into the complement of $\nu$ by such $\cdot 11\cdot\to \cdot 2\cdot$ moves. This gives the reduction
$1^m\nu \to 1^{m'}\nu^c\nu$ where $m' \ge 3$.

The $\nu$ sequence and its complement can be glued up in an obvious way, folding up the segment between
the last letter of $\nu^c$ and the first letter of $\nu$ and producing an interior orbifold point of order
$2$. This amounts to the reduction $1^{m'}\nu^c\nu \to 1^{m''}2$ where $m'' \ge 1$. After a finite
sequence of $\cdot 1121\cdot \to \cdot 2\cdot$ moves, we reduce to one of the cases $12$, $112$ or $1112$.

\begin{figure}[htpb]
\labellist
\small\hair 2pt
\pinlabel $\text{fold}$ at -10 41
\pinlabel $1$ at 40 12
\pinlabel $2$ at 40 69
\pinlabel $2$ at 104 64
\pinlabel $1$ at 104 17
\pinlabel $1$ at 146 45
\pinlabel $1$ at 200 12
\pinlabel $2$ at 200 69
\pinlabel $1$ at 174 45
\pinlabel $1$ at 226 45
\endlabellist
\centering
\includegraphics[scale=1]{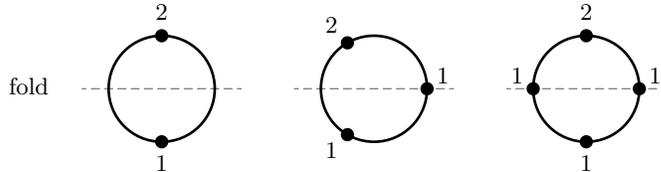}
\caption{Gluing up $12$, $112$ and $1112$ by folding}\label{12_112_1112}
\end{figure}

Folding up each edge of the $12$ circle glues up the boundary completely, producing two interior
orbifold points of order $2$. Folding up the $\cdot 12\cdot$ edge and identifying the other pair of edges in
the $112$ circle glues up the boundary completely, producing two interior orbifold points or orders
$2$ and $3$. Identifying the $\cdot 12\cdot$ and $\cdot 21\cdot$ edges with the succeeding pair of $\cdot 11\cdot$ edges of the
$1112$ circle glues up the boundary completely, producing two interior orbifold points of order $3$.
See Figure~\ref{12_112_1112}.

In every case therefore $\partial_W S'$ can be completely glued up, and the proof of 
Theorem~\ref{stability_theorem} is complete.

\begin{remark}
The proof generalizes in an obvious way to formal sums. Let $\gamma_n$ denote the geodesic corresponding
to the conjugacy class $R^nX$ for some fixed string $X$. Let $\delta_1,\cdots,\delta_m$ be any finite
collection of geodesics in $\M$. Then for sufficiently large $n$, the $1$-manifold $\gamma_n \cup_i \delta_i$
virtually bounds an immersed surface in $\M$.
\end{remark}

\section{Generalizations}\label{generalize_section}

In this section we discuss some generalizations of our main theorem.

\subsection{$(2,p,\infty)$-orbifold}

The purpose of this section is to give a proof of the following theorem which generalizes
Theorem~\ref{stability_theorem} to hyperbolic $(2,p,\infty)$-orbifolds for any $p\ge 3$.

\begin{theorem}[$(2,p,\infty)$-stability]\label{generalized_stability}
Let $\M_p$ denote the hyperbolic $(2,p,\infty)$-orbifold for $p\ge 3$, and let $G_p$ be its fundamental
group. Let $R \in G_p$ be the element corresponding to a positive loop around the puncture. Let
$v \in G_p$ be arbitrary. Then for all sufficiently large $n$, the geodesic in $\M_p$ corresponding
to the stabilization $R^nv$ virtually bounds an immersed surface in $\M_p$ (providing this element
is hyperbolic).
\end{theorem}

\begin{remark}
For any $v$, the elements $R^nv$ are eventually hyperbolic unless $v$ is a power of $R$.
\end{remark}

\begin{proof}
The proof is very similar to the proof of Theorem~\ref{stability_theorem}, and for the sake of
brevity we use language and notation as in \S~\ref{proof_section} by analogy.

There is an embedded geodesic segment $\sigma_p$ running between the orbifold points of orders $2$ and
$p$ in $\M_p$, covered by an infinite $p$-valent tree $\til{\sigma}_p$ in $\H^2$. Let $W_p$ be
the closure of the complementary component of $\H^2 -\til{\sigma}_p$ stabilized by a translation.
Then $\partial W_p$ is a sequence of circular arcs meeting at an angle of $2\pi/p$.

A conjugacy class $R^nv$ represented by a geodesic $\gamma$ is decomposed into $\alpha_i$ and
$\beta_i$ arcs by $\sigma_p$, where we use $\alpha_i$ to denote the arcs whose lifts to $W_p$ move
to the left, and the $\beta_i$ to denote the arcs whose lifts to $W_p$ move to the right. For
sufficiently large $n$, there is one arc $\alpha_1$ with a lift $\til{\alpha}_1$
which intersects segments of $\partial W_p$ approximately $n$ apart, whereas the combinatorial
types of the $\alpha_i$ (for $i>1$) and $\beta_i$ are eventually constant.

As in \S~\ref{lifts_surfaces_subsection} we choose disjoint lifts $\til{\alpha}_i$, $\til{\beta}_i$
with the $\til{\beta}_i$ under $\til{\alpha}_1$ and with five segments of $\partial W_p$ between
successive $\til{\beta}_i$. These lifts cobound surfaces $T$ and $S_i$ which can be glued up to
produce $S'$ with $\partial_W S'$ a collection of circles with vertices labeled by numbers from
$1$ to $(p-1)$. As in \S~\ref{lifts_surfaces_subsection}, we are guaranteed that each circle contains 
a $\cdot 11\cdot$, and one circle contains a big $1$ sequence (which may be assumed to be as long as we like,
by making $n$ large).

The $1$-handle move still makes sense, so after performing a sequence of such moves we can reduce
to the case of a single circle with a big $1$ sequence; i.e. we have reduced to the case of
$1^m\nu$ for $\nu$ an arbitrary (but fixed independent of all large $n$) sequence, and $m$ as big
as we like. Folding a $\cdot 1k\cdot$ segment where $k<p-1$ produces an interior orbifold point of order $2$, and
a single $(k+1)$ vertex. So we can turn a sequence of at most $(p-1)|\nu|$ consecutive $1$'s into the
complement $\nu^c$ and reduce $1^m\nu \to 1^{m'}\nu^c\nu$. Folding $\nu^c$ into $\nu$ reduces to
$1^{m''}2$ as before.

We must be a bit careful with the endgame depending on the value of $m''$ mod $p$. We would like to reduce
to a $1^{p-k}k$ circle for some $k$, since
$1^{p-k}k \to 1^{p-k-1}(k+1) \to \cdots \to 1(p-1)$ by successive folds, and then $1(p-1)$ can be completely
folded up, producing two orbifold $2$ points (much as we completely folded up the $12$ circle in the
case $p=3$).

Fortunately this reduction can be accomplished if $m''$ is sufficiently big (whatever its residue mod $p$).
Folding an edge gives $\cdot 11\cdot \to \cdot 2\cdot$, but folding at a vertex gives $\cdot 111\cdot \to \cdot 2\cdot$,
in either case producing an interior orbifold point of orders $2$ and $p$ respectively. By
judicious application of some number of these moves (together with folds of the kind 
$\cdot a b\cdot \to \cdot (a+b) \cdot$ if $a+b<p$ or $\cdot 1ab1\cdot \to \cdot 2\cdot$ if $a+b=p$)
we can reduce $1^{m''}2$ for any sufficiently large $m''$ to $1^{p-k}k$ for some $k$, and thence
glue up completely as above. This completes the proof of the theorem.
\end{proof}

It seems very plausible that some generalization of our methods should prove an analogous theorem
for $(p,q,\infty)$ orbifolds with $p,q>2$, or even arbitrary noncompact hyperbolic orbifolds 
with underlying topological space a disk, but the combinatorial endgame becomes progressively more
complicated, and we have not pursued this.

\subsection{Combinatorics of $12$ circles}

Returning to the case of $p=3$, we describe a slightly different method for producing immersed surfaces.
Suppose we have a word $X=R^{a_1}L^{b_1}R^{a_2}\cdots L^{b_n}$ for which the $b_i$ can be partitioned
into subsets $B_i = \lbrace b_{i,1}, b_{i,2},\cdots \rbrace$ (possibly empty) with the property that
$a_i \gg \sum_j b_{i,j}$. In this case, we can choose lifts $\til{\alpha}_i$ and $\til{\beta}_j$
such that for each $i$, all the $\til{\beta}_{i,j}$ are contained ``under'' the single arc $\til{\alpha}_i$.

In this case, we can still produce a surface $S'$ for which $\partial_W S'$ is a collection of 
circles with $1$ and $2$ vertices. We are thus naturally led to the question: which $12$ circles
can be glued up completely? This seems like a hard combinatorial question; nevertheless, we describe
some interesting necessary conditions and sufficient conditions (though we don't know a simple condition
which is both necessary and sufficient).

\begin{example}
Each $2$ vertex must be glued to some unique $1$ vertex, therefore the total number of $2$'s must be
no more than the total number of $1$'s. So, for example, the circles $2$ and $221$ can't be glued up.
\end{example}

\begin{example}
The length of a consecutive sequence of $1$'s can never be increased. So consecutive sequences of
$2$'s must be associated to disjoint consecutive sequences of $1$'s of at least the same length. So,
for example, the circle $22211211211$ can't be glued up, even though it has more $1$'s than $2$'s.
\end{example}

\begin{example}
A circle with alternating $1$'s and $2$'s can be completely glued up. Any circle
of the form $21^{c_1}21^{c_2}\cdots 21^{c_m}$ where each $c_i\ge 7$ can be completely glued up,
since we can use the reductions $\cdot 111\cdot \to \cdot 2 \cdot$ or $\cdot 11 \cdot \to \cdot 2 \cdot$
to replace $\cdot 1^c\cdot $ by $\cdot 121\cdots 21\cdot$ whenever $c\ge 7$. Hence in particular,
$\rot$ is extremal for any conjugacy class of the form $R^{a_1}L^{b_1}R^{a_2}L\cdots R^{a_n}L^{b_n}$ 
whenever the $b_i$ can be partitioned into subsets $B_i = \lbrace b_{i,1}, b_{i,2},\cdots\rbrace$
(possibly empty) with the property that $a_i \ge 10+\sum_j (b_{i,j}+10)$ for all $i$.
\end{example}

Given two finite sets of positive numbers $\lbrace a_i\rbrace$ and $\lbrace b_i \rbrace$
(possibly with multiplicity), the problem of partitioning
the $b_i$ into subsets $B_i = \lbrace b_{i,1}, b_{i,2},\cdots \rbrace$ (possibly empty) with
the property that $a_i \ge \sum_j b_{i,j}$ is familiar in computer science, where the $b_i$ denote
the lengths of a family of files, and the $a_i$ denote the lengths of a family of empty consecutive
blocks in memory. See e.g. \cite{Knuth}, \S~2.2 and \S~2.5 for a discussion. The performance of
dynamic memory allocation algorithms is very well studied, with respect to many different kinds of 
statistical distributions for the numbers $\lbrace a_i\rbrace$ and $\lbrace b_i \rbrace$, 
and it would be intriguing to pursue this connection further.

\begin{example}
In this paper we have discussed various sufficient conditions on the exponents 
$a_i$ and $b_i$ for a word $\gamma$ to virtually bound an immersed surface. However, 
these conditions have only depended on the {\em sets} of values $\lbrace a_i \rbrace$ and
$\lbrace b_i \rbrace$, and not their {\em order}. A complete understanding must necessarily
take this order into account. For example,
$$\scl(R^3LRLR^2LRL^2) = 1/6 = \rot(R^3LRLR^2LRL^2)/2$$ 
whereas
$$\scl(R^3LRLRLR^2L^2) = 4/15 > 1/6 = \rot(R^3LRLRLR^2L^2)/2$$
\end{example}

\section{Experimental results}\label{experimental_section}

In this section we describe the results of some computer experiments, comparing the functions
$\scl(\cdot)$ and $\rot(\cdot)/2$ in general. The function $\rot$ can be computed by an exponent sum
for a conjugacy class expressed as a product of $R$'s and $L$'s, and $\scl$ can be computed using the algorithms
described in \cite{Calegari_rational} (implemented on the program {\tt scallop}, available from
\cite{Calegari_scallop}) and \cite{Calegari_scl}, \S~4.2.5. 

\subsection{Distribution of $n(X)$}

For each word $X$ in $L$ and $R$, define $n(X)$ to be the smallest negative number such that
$\rot$ is extremal for $L^{-n}X$ (if one exists), or the smallest non-negative number such that
$\rot$ is extremal for $R^nX$ otherwise.

\begin{figure}[htpb]
\labellist
\small\hair 2pt
\pinlabel $-10$ at 20 0
\pinlabel $-5$ at 165 0
\pinlabel $0$ at 310 0
\pinlabel $5$ at 455 0
\pinlabel $\ge 10$ at 600 0
\endlabellist
\centering
\includegraphics[scale=0.5]{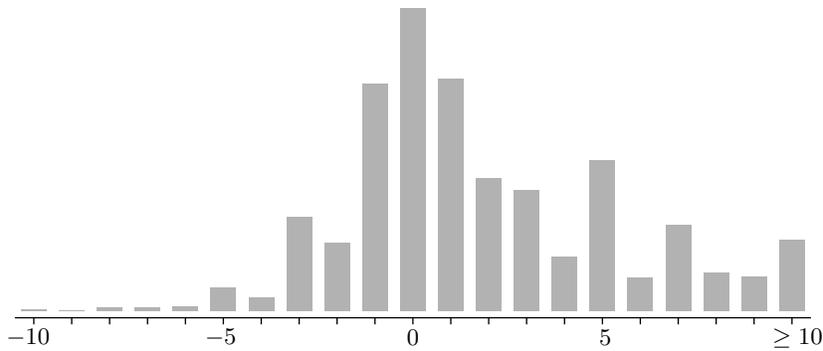}
\caption{Histogram of $\freq(n(X))$ for all 1024 words $X$ of length $10$}\label{histogram_10}
\end{figure}

Figure~\ref{histogram_10} shows a histogram of the frequency distribution of $n(X)$, for all
words $X$ of length $10$.

It is a fact (\cite{Calegari_Maher}) that in an arbitrary word-hyperbolic group, most
rationally null-homologous words of length $n$ have $\scl \sim n/\log{n}$. On the other hand,
$\rot$ is an example of a {\em bicombable quasimorphism}, and therefore by \cite{Calegari_Fujiwara},
the distribution of values on words of length $n$ satisfies a central limit theorem; in 
particular, one has $|\rot| \sim \sqrt{n}$ for most words of length $n$. 
It follows that $n(X)$ is at least of size $n/\log{n}$ for most words $X$ of length $n$,
at least when $n$ is large.

\subsection{Stuttering}

One might imagine from the discussion above that if $\rot$ is extremal for $R^nX$, then
it must also be extremal for $R^mX$ for all $m>n$; however, this is not the case. We call
this phenomenon {\em stuttering}.

\begin{example}
The quasimorphism $\rot$ is extremal for $R^3LRL^2$ but not for $R^4LRL^2$. It is extremal
for $R^2LRL^2RL$ but not for $R^3LRL^2RL$ or $R^4LRL^2RL$. It is extremal
for $RLR^2L^2RLRL^2R^2L$ but not for $R^iLR^2L^2RLRL^2R^2L$ for $1<i<5$. We refer to these examples
colloquially as {\em stuttering sequences} of length $1$, $2$ and $3$ respectively. We do not
know of any examples of stuttering sequences of length $>3$, but do not know any reason why
such examples should not exist.
\end{example}

\section{Acknowledgments}
Danny Calegari was supported by NSF grant DMS 0707130. We would like to thank Benson Farb and
Eric Rains for some useful conversations about this material.


\begin{thebibliography}{99}
\bibitem{Atiyah}
  M. Atiyah,
  \emph{The logarithm of the Dedekind $\eta$-function},
  Math. Ann {\bf 278}, 1-4 (1987), 335--380
\bibitem{Bavard}
  C. Bavard,
  \emph{Longeur stable des commutateurs},
  L'Enseign. Math. {\bf 37} (1991), 109--150
\bibitem{Burger_Iozzi_Wienhard}
  M. Burger, A. Iozzi and A. Wienhard,
  \emph{Surface group representations with maximal Toledo invariant},
  C. R. Math. Acad. Sci. Paris {\bf 336} (2003), no. 5, 387--390
\bibitem{Calegari_rational}
  D. Calegari,
  \emph{Stable commutator length is rational in free groups},
  Jour. AMS {\bf 22} (2009), no. 4, 941--961
\bibitem{Calegari_faces}
  D. Calegari,
  \emph{Faces of the scl norm ball},
  Geom. Topol. {\bf 13} (2009), no. 3, 1313--1336
\bibitem{Calegari_scl}
  D. Calegari,
  \emph{scl},
  MSJ Memoirs, {\bf 20}. Mathematical Society of Japan, Tokyo, 2009
\bibitem{Calegari_scallop}
  D. Calegari,
  {\tt scallop}, computer program available from the author's webpage, and from {\tt computop.org}
\bibitem{Calegari_Fujiwara}
  D. Calegari and K. Fujiwara,
  \emph{Combable functions, quasimorphisms and the central limit theorem},
  Ergodic Theory Dynam. Systems, to appear
\bibitem{Calegari_Maher}
  D. Calegari and J. Maher,
  in preparation
\bibitem{Ghys}
  \'E. Ghys,
  \emph{Knots and dynamics},
  Internat. Congress Math. Vol. 1, 247--277 Eur. Math. Soc., Z\"urich 2007
\bibitem{Hirzebruch_Zagier}
  F. Hirzebruch and D. Zagier,
  \emph{The Atiyah-Singer theorem and elementary number theory},
  Math. Lect. Ser. 3, Publish or Perish, Boston 1974
\bibitem{Kirby_Melvin}
  R. Kirby and P. Melvin,
  \emph{Dedekind sums, $\mu$-invariants and the signature cocycle},
  Math. Ann. {\bf 299} (1994), 231--267
\bibitem{Knuth}
  D. Knuth,
  \emph{The art of computer programming. Vol. 1: Fundamental algorithms},
  Second Printing Addison-Wesley Publishing Co., Reading, Mass. 1969
\bibitem{Rademacher}
  H. Rademacher and E. Grosswald,
  \emph{Dedekind sums},
  Carus Mathematical Monographs 16, Math. Assoc. Amer., Washington DC 1972
\bibitem{Thurston}
  W. Thurston,
  \emph{A norm for the homology of $3$-manifolds},
  Mem. Amer. Math. Soc. {\bf 59} (1986), no. 339, i--vi and 99--130
\bibitem{Wilton}
  H. Wilton,
  private communication
\end{thebibliography}
\end{document}